# DYNAMIC TREE ALGORITHMS

By Hanène Mohamed and Philippe Robert

*Universite de Versailles-saint-Quentin-en-Yvelines and INRIA*

*This paper is dedicated to Philippe Flajolet on the occasion of his 60th birthday*

In this paper, a general tree algorithm processing a random flow of arrivals is analyzed. Capetanakis–Tsybakov–Mikhailov's protocol in the context of communication networks with random access is an example of such an algorithm. In computer science, this corresponds to a trie structure with a dynamic input. Mathematically, it is related to a stopped branching process with exogenous arrivals (immigration). Under quite general assumptions on the distribution of the number of arrivals and on the branching procedure, it is shown that there exists a *positive* constant $\lambda_c$ so that if the arrival rate is smaller than $\lambda_c$, then the algorithm is stable under the flow of requests, that is, that the total size of an associated tree is integrable. At the same time, a gap in the earlier proofs of stability in the literature is fixed. When the arrivals are Poisson, an explicit characterization of $\lambda_c$ is given. Under the stability condition, the asymptotic behavior of the average size of a tree starting with a large number of individuals is analyzed. The results are obtained with the help of a probabilistic rewriting of the functional equations describing the dynamics of the system. The proofs use extensively this stochastic background throughout the paper. In this analysis, two basic limit theorems play a key role: the renewal theorem and the convergence to equilibrium of an auto-regressive process with a moving average.

**1. Introduction.** This paper investigates probabilistic algorithms which decompose recursively a given set of elements (also referred to as items) into random subsets until all subsets have a cardinality less than some fixed number $D$. The dynamic aspect of the algorithms analyzed here is that new elements are added to the subsets created during each decomposition.

The general procedure is as follows: if the cardinality of the set is strictly less than $D > 0$, the process is stopped. Otherwise, the set is split into several subsets and each subset receives a random number of new elements. The algorithm is then recursively applied to each of these subsets. A tree is naturally associated with this algorithm: the root node having the initial items; the subsequent nodes containing the corresponding subsets, and so on. At the end of this process, one ends up with a tree whose leaves contain less than $D$ items and all internal nodes contain more than $D$ items. Such an algorithm can be seen as a state dependent branching process which dies out whenever a termination condition is satisfied. When there are no arrivals of new elements, the algorithm is called *static*. Static tree algorithms are of fundamental importance as a generic class. They are used in computer science where the corresponding data structure is called a *trie*, and in many other areas such as communication protocols and distributed systems. See Mohamed and Robert [14] for a general overview of static tree algorithms.

The extension analyzed here, with the Introduction of new elements, consists of introducing immigration to the language of branching processes. This situation is quite natural in the context of communication protocols where new requests (immigration) arrive continuously at the communication node. On the use of these algorithms in the context of communication networks, see the surveys by Berger [2] and Massey [12]. On the mathematical side, as it will be seen, the analysis of these algorithms turns out to be more delicate.

TREE ALGORITHM WITH IMMIGRATION $\mathcal{S}(n)$.

---

– TERMINATION CONDITION.
   If $n < D$          $\longrightarrow$ STOP.
– TREE STRUCTURE.
   If $n \geq D$, randomly divide $n$ into $n_1, \ldots, n_G$, with $n_1 + \cdots + n_G = n$.
                     $\longrightarrow$ APPLY $\mathcal{S}(n_1 + A_1), \mathcal{S}(n_2 + A_2), \ldots, \mathcal{S}(n_G + A_G)$
   where $(A_i)$ are i.i.d. random variables.

---

For the static algorithm, when there are no new arrivals, provided that the decomposition mechanism is not degenerated, it is easily seen that the associate tree is almost surely finite, in fact, that the total number of its nodes is integrable. Mohamed and Robert [14] investigates this case.

*Finiteness of the associated tree and law of large numbers.* When there is a set of new items arriving with every time unit, it may happen that the algorithm does not terminate with probability 1, that is, that the associated tree is infinite. In this case, the algorithm cannot cope with the flow of arriving requests. This nontrivial phenomenon is analyzed in this paper. Furthermore, in the case where the process terminates almost surely and that there are $n$ items at the root, another problem is to describe the asymptotic behavior of the average size of the tree as $n$ gets large. As we will see, this is a quite challenging question.

*Tree recurrences.* In this setting the main quantity of interest $R_n^A$ is the total number of nodes of the associated tree when the algorithm starts with $n$ items. The superscript "$A$" of $R_n^A$ refers to the common distribution of the i.i.d. sequence $(A_i)$. In this case, one gets naturally the recursive relation

$$(1) \qquad R_n^A \stackrel{\text{dist}}{=} 1 + (R_{1,n_1+A_1}^A + R_{2,n_2+A_2}^A + \cdots + R_{G,n_G+A_G}^A)\mathbf{1}_{\{n \geq D\}},$$

where, for $i \geq 1$, $(R_{i,n}^A, n \geq 0)$ are independent random variables with the same distribution as $(R_n^A, n \geq 0)$. The total size of the tree is 1 plus the size of all of sub-trees of level 1.

*A Markovian representation.* This algorithm can also be represented by a Markov chain $(\underline{L}_t)$ on the set $\mathcal{S} = \bigcup_{\ell \geq 0} \mathbb{N}^\ell$ of finite sequences on $\mathbb{N}$; its transitions are described as follows: if $(\underline{L}_0) = (l_0, l_1, l_2, \dots)$

$$(2) \qquad (\underline{L}_1) = \begin{cases} (l_1 + A_1, l_2, \dots, l_n, \dots) \text{ if } l_0 < D, & /\text{SSHIFT}/, \\ (n_1 + A_1, n_2, \dots, n_G, l_1, l_2, l_3, \dots) \text{ if } l_0 \geq D, & /\text{SPLIT}/, \end{cases}$$

if the integer $l_0$ is decomposed into $l_0 = n_1 + n_2 + \cdots + n_G$ by the splitting procedure (see the precise description below).

Note that if (1) and (2) are representations of the dynamics of tree algorithm; they differ in the following way. Equation (2) gives the state of the system just after one time unit, that is, the number of individuals with counter $\ell$, $\ell \geq 0$, including the $A_1$ new items. Equation (1) is a branching representation for the size of the "final" tree with its $G$ children. For $1 \leq i \leq G$, $A_i$ is the number of new arrivals at the beginning of the time unit when the value of the counter of the $i$th subgroup (with cardinality $n_i$) is 0.

For $n \geq 1$, if the initial state is $(\underline{L}_0 = (n, 0, \dots, 0, \dots))$, then it is not difficult to check that equation (1)

$$R_n^A \stackrel{\text{dist}}{=} \inf\{t \geq 1 : (\underline{L}_t = (0, \dots, 0, \dots))\}.$$

The variable $R_n^A$ can also be viewed as the hitting time of the empty state by the Markov chain $(\underline{L}_t)$. The ergodicity of the Markov chain $(\underline{L}_t)$ is, therefore, equivalent to the fact that the variable $R_{A_1}^A$ is integrable. Because of the description by a sequence (a *stack* in the language of computer science), the tree algorithm is, sometimes, also called *stack algorithm*.

Provided that the variables $(R_n^A)$ are integrable, the Poisson transform $\phi(x)$ of the sequence $(\mathbb{E}(R_n^A))$ is defined as

$$(3) \qquad \phi(x) \stackrel{\text{def}}{=} \sum_{n \geq 0} \mathbb{E}(R_n^A)\frac{x^n}{n!}e^{-x}.$$

In the case where the arrivals have a Poisson distribution with parameter $\lambda$, the ergodicity of $(\underline{L}_n)$ implies that the Poisson transform of the sequence $(\mathbb{E}(R_n^A))$ is well defined at $x = \lambda$.

*Iteration of noncommutative functional operators.* Mathematically, these tree algorithms are quite challenging. By using an iterating scheme, a family of simple functional operators $(P_v^\lambda, v \in (0,1))$ play a central role in the analysis; if $f \colon \mathbb{R}_+ \to \mathbb{R}_+$ is a continuous function, they are defined by

$$P_v^\lambda(f)(x) = \frac{1}{v} f(xv + \lambda).$$

The static algorithm corresponds to the case when $\lambda = 0$.

By taking the expected value of (1) and by iterating the functional obtained, it turns out that the Poisson transform can be expressed by the following equation:

$$(4) \qquad \sum_{n=1}^{+\infty} \int_{[0,1]^n} P_{v_1}^\lambda \circ P_{v_2}^\lambda \circ \cdots \circ P_{v_n}^\lambda (f)(x) \bigotimes_{i=1}^n \mathcal{W}(dv_i), \qquad x \geq 0,$$

for a convenient function $f$ depending on some unknown constants and where $\mathcal{W}$ is some probability distribution on $(0,1)$.

Note that the operators $(P_v^\lambda, v \in [0,1])$ commute only when $\lambda = 0$. This complicates significantly the analysis of the asymptotic behavior of expression (4) as $x$ goes to infinity. In the static case, one has that

$$P_{v_1}^0 \circ P_{v_2}^0 \circ \cdots \circ P_{v_n}^0 = P_{v_1 v_2 \cdots v_n}^0,$$

which gives a multiplicative representation of expression (4) which is analyzed by using Mellin transform methods in an analytical context (see Flajolet, Gourdon and Dumas [10]) or by using random walks methods with a probabilistic approach (see Mohamed and Robert [14]). When $\lambda > 0$, and such a multiplicative formulation is not available, these methods have, therefore, to be adapted. It turns out that such a generalization is not straightforward.

EXAMPLE (The binary tree). The splitting mechanism is binary, and the branching number $G$ is deterministic and equal to 2, $G \equiv 2$, and with deterministic weights, $V_{11} \equiv p$ and $V_{12} \equiv q = 1 - p$ with $p \in (0,1)$. The variables $(A_i)$ are assumed to be Poisson with parameter $\lambda$. The upper index $A$ in $(R_n^A)$ is replaced by $\lambda$ in this case. Provided that the variables are integrable, and if $\alpha_n^\lambda = \mathbb{E}(R_n^\lambda)$, the integration of (1) gives the identities $\alpha_0^\lambda = \alpha_1^\lambda = 1$, and, for $n \geq 2$,

$$(5) \qquad \alpha_n^\lambda = 1 + \sum_{i=0}^n \binom{n}{i} p^i q^{n-i} \sum_{k,\ell \geq 0} \frac{\lambda^k}{k!} e^{-\lambda} \frac{\lambda^\ell}{\ell!} e^{-\lambda} (\alpha_{i+k}^\lambda + \alpha_{n-i+\ell}^\lambda).$$

In this case it is not difficult to check that the corresponding Poisson transform $\phi$ satisfies the equation

$$(6) \qquad \phi(x) = \phi(px + \lambda) + \phi(qx + \lambda) + h(x),$$

where $x \to h(x)$ is some fixed function with a specific form and unknown coefficients. This functional relation can be rewritten as

$$\phi(x) = \int_0^1 P_w^\lambda(\phi)(x)\mathcal{W}(dw) + h(x),$$

where $\mathcal{W} = p\delta_p + q\delta_q$, where $\delta_x$ is the Dirac measure at $x \in \mathbb{R}$. A (formal) iteration of this function gives representation (4) for $\phi$.

*Literature.* These problems have been analyzed in several ways in the past. Motivated by the design of stable communication protocols, Tsybakov and Mikhaïlov [15] and Capetanakis [6] did the early studies in this domain (and, at the same time, designed the algorithms in the context of distributed systems) (see also Tsybakov and Vvedenskaya [16]).

At the end of the 1980s, Flajolet and his co-authors, in a series of interesting papers [8, 9, 13], have obtained rigorous asymptotics for solutions of the type of equations as seen in (5), in several cases. In the first of the papers [8], recurrence (5) for the binary tree is investigated. It is shown that for $\lambda$ smaller than some threshold, there is a unique sequence $(\alpha_n^\lambda)$ of real numbers which is the solution of this recurrence. In this paper, a sophisticated asymptotic analysis of the sequence $(\alpha_n^\lambda)$ is presented. Basically, it is shown that the sequence $(\alpha_n^\lambda)$ grows linearly with respect to $n$ and its rate is, in some cases, a function with small fluctuations. It is conducted in three steps:

(1) by iteration of this equation, express the Poisson transform $\phi$ of $(\alpha_n^\lambda)$ under the form (4),

$$\tag{7} \phi(x) = \sum_{n \geq 0} \sum_{i \in \{1,2\}^n} h(\sigma_{i_1} \circ \sigma_{i_2} \circ \cdots \circ \sigma_{i_n}(x)),$$

where $\sigma_1(x) = px + \lambda$ and $\sigma_2(x) = qx + \lambda$;
(2) obtain the asymptotics of the Poisson transform $\phi(x)$ as $x$ goes to infinity;
(3) prove that $(\phi(x))$ and $(\alpha_n^\lambda)$ have the same behavior as $x$ (resp. $n$) goes to infinity.

The main part of the analysis is devoted to step (2) where, via several estimates of series and contour integrals with arguments from complex analysis, the authors can identify in the series (7) the main contributing terms when $x$ goes to infinity. Given the complexity of this analysis for this example of the binary tree, an extension of these methods to a more general branching mechanism seems to be more than challenging. The case when $\mathcal{W}$ is the uniform distribution on $Q$ points is discussed in Section IV of Mathys and Flajolet [13]. Note that if $\mathcal{W}$ has some Lebesgue component, and a representation of $\phi$ as a series similar to (7) is no longer available, one has to go back to the general representation (4).

In addition to extending these results to a quite general branching scheme, the purpose of this paper is also to develop probabilistic methods to analyze additive functionals of various tree structures. This program has been initiated in Mohamed and Robert [14] in the context of static tree algorithms. We nevertheless believe that, because of the intricacies of its associated equations, the tree algorithm with arrivals provides a real significant test for this approach. It turns out that the method proposed in this paper has some advantage in that it can handle more easily and in a more general setting the complexity of the underlying noncommutative iterating scheme of this algorithm.

Recent works deal with some aspects of these fundamental algorithms, see Boxma, Denteneer and Resing [5], Janssen and de Jong [11] and Van Velthoven, Van Houdt and Blondia [17], for example. For surveys on the communication protocols in random access networks, see Berger [2], Massey [12] and Ephremides and Hajek [7].

*Contributions of the paper.*   The main objective of the paper is to present a probabilistic approach to these problems that can tackle, with a limited technical complexity, quite general models of tree algorithms with immigration.

For the model considered in this paper, (6) of the Poisson transform becomes

$$(8) \qquad \phi(x) = \int_0^1 \frac{\phi(wx + \lambda)}{w} \mathcal{W}(dw) + h(x)$$

for some probability distribution $\mathcal{W}$ on $(0,1)$ describing the branching mechanism of the splitting algorithm and some function $h$. The example of the binary tree corresponds to the case where $\mathcal{W}$ is carried by two points $p$ and $q$, as noted before. When the measure $\mathcal{W}$ is carried by $Q$ points of $(0,1)$, the equivalent expression for the series (7) involves the various products of $Q$ functions $(\sigma_m, 1 \leq m \leq Q)$.

*A. Stability of tree algorithms.*   The stability property of the tree algorithm, that is, the fact that the tree is almost surely finite, is a crucial issue for communication networks. Assuming Poisson arrivals with parameter $\lambda$, it amounts to the existence of some threshold $\lambda_0 > 0$ such that if the arrival rate $\lambda$ is strictly less than $\lambda_0$, then the associated Markov process describing the tree algorithm (see above) is ergodic. During the 1970s and 1980s the design of stable protocols and the mathematical proof of their stability has been a very active research domain. Recently, because of the emergence of wireless and mobile networks, there is a renewed interest in these models. The first protocols, Aloha and Ethernet, turned out to be unstable when there is an infinite number of possible sources, that is, for these algorithms, the number of requests waiting for transmission goes to infinity in distribution for any arrival rate (see Aldous [1]). The tree algorithm corresponding to the example of the binary tree with $p = q = 1/2$ is the first such protocol in a really distributed system where the stability region is nonempty (see Massey [12] and Bertsekas and Gallager [4]). In later papers, the tree algorithm has been improved in order to have a larger stability region.

*A gap in the proof of previous stability results.* The stability results obtained up to now are under the assumption of Poisson arrivals. The proofs known to us rely on the analysis of of the type as in equation as in (5) for the sequence $(\mathbb{E}(R_n^\lambda))$; it is shown that there exists some $\lambda_0 > 0$ such that for $\lambda < \lambda_0$, there is a unique finite solution $(\alpha_n)$ and from there it is concluded that the corresponding Markov chain is ergodic. The problem here is that this analysis shows only that, for $\lambda < \lambda_0$, there exists a unique sequence $(\alpha_n)$ of finite real numbers satisfying relation (5). The sequence $(\beta_n) = (1, 1, +\infty, \ldots, +\infty, \ldots)$ also satisfies relation (5). At this point, without an additional argument, it cannot be concluded that, for $\lambda < \lambda_0$, $(\mathbb{E}(R_n^\lambda))$ is indeed $(\alpha_n)$, and not $(\beta_n)$. To make the identification with $(\alpha_n)$, it has to be proved that the random variables $R_n^\lambda$, $n \in \mathbb{N}$, are indeed integrable, but this is precisely the final result. Strictly speaking, the previous results have only shown that the system is unstable whenever $\lambda \geq \lambda_0$.

This gap is fixed in this paper. The key ingredient to relate recurrence in relation (5) and the sequence $(R_n^\lambda)$ is a perturbation result of the static case, that is, when $\lambda = 0$. It is also shown that, under general assumptions on the distribution of the inputs $(A_i)$ and on the branching mechanism, the tree algorithm is stable for a sufficiently small arrival rate. As far as we know, this is the first stability result for tree algorithms with non-Poissonnian arrivals.

*B. Analysis of general tree recurrences.* The second part of the paper investigates the asymptotic behavior of the sequence $(\mathbb{E}(R_n^A)/n)$ where $(R_n^A)$ is a solution of the tree recurrence (1) under the condition that $A$ is a Poisson random variable with a parameter $\lambda$ less than some constant.

Some of the ingredients of the analysis of static algorithms ($\lambda = 0$) of Mohamed and Robert [14] are used. The situation is nevertheless completely different when $\lambda > 0$. As mentioned above, the noncommutativity of the operators $P_v$ is a major issue. An additional important difficulty is the fact that, contrary to the case $\lambda = 0$, the function $h$ of (8) is unknown, and $D$ coefficients have to be determined.

To study these recurrences, the approach of the paper consists of expressing series (7) for the binary tree or (4) in the general case as the expected value of some random variable depending on some auto-regressive process with moving average $(X_n)$ defined by $X_0 = x$ and

$$(9) \qquad X_{n+1} = W_n X_n + 1, \qquad n \geq 0,$$

where $(W_n)$ is an i.i.d. sequence whose common distribution is $\mathcal{W}$.

Key limit theorems are used to derive the asymptotic behavior of the sequence $(\mathbb{E}(R_n^\lambda))$ such as the renewal theorem and the convergence of $(X_n)$ to its stationary distribution. In the particular case of the binary tree, they avoid the use of estimations of Fayolle, Flajolet and Hofri [8] which are necessary to get the significant terms of the series (7) in the asymptotic expansion. Roughly speaking, with these two limit theorems, the probabilist "knows" what are the most likely trajectories of the compositions of $\sigma_1$ and $\sigma_2$. This approach simplifies *significantly* the asymptotic analysis of the algorithm. Another key step is to identify the $D$ unknown coefficients of function $h$; they are expressed as a functional of the auto-regressive process, of its invariant distribution in particular.

*Outline of the paper.* The paper is organized as follows: Section 2 shows that, under a quite general assumption of arrivals, the stability region is not empty. In Section 3, by denoting $R_n^\lambda$, the size of the tree when $n \geq 0$ elements are at the root, and under the assumption of Poisson arrivals, a probabilistic representation of the Poisson transform of the sequence $(\mathbb{E}(R_n^\lambda))$ is established, and an auto-regressive process with a moving average is introduced. Section 4 establishes the main stability result (Theorem 9) for general tree algorithm Poisson arrivals. Section 5 investigates the delicate asymptotics of the sequence $(\mathbb{E}(R_n^\lambda)/n)$, and Theorem 12 summarizes the results obtained.

**2. Existence of a nonempty stable region.** In this section it is proved that, if the arrival rate is sufficiently small, then the tree obtained with the algorithm is almost surely finite; its size is in fact integrable.

*Formulation of the problem.* The algorithm starts with a set of $n$ items. If $n < D$, then it stops. Otherwise, this set is randomly split into $G$ subsets where $G$ is some random variable. Now, conditionally on the event $\{G = \ell\}$, for $1 \leq i \leq \ell$, each of the $n$ items is sent into the $i$th subset with probability $V_{i,\ell}$ where $V_\ell = (V_{i,\ell}; 1 \leq i \leq \ell)$ is a random probability vector on $\{1, \ldots, \ell\}$. The quantity $V_{i,\ell}$ is the weight on the $i$th edge of the splitting structure. Additionally, a vector $(A_1, \ldots, A_\ell)$ of independent random variables with the same distribution as some random variable $A$ is given, and $A_i$ new items are added to the $i$th subset. If $n_i$ is the cardinality of the $i$th subset, then, conditionally on the event $\{G = \ell\}$ and on the random variables $V_{1,\ell}, V_{2,\ell}, \ldots, V_{\ell,\ell}$, the distribution of the vector $(n_1, \ldots, n_\ell)$ is multinomial with parameter $n$ and $(V_{1,\ell}, V_{2,\ell}, \ldots, V_{\ell,\ell})$. If the $i$th subset, $1 \leq i \leq n$, is such that $n_i + A_i < D$, the algorithm stops for this subset. Otherwise, it is applied to the $i$th subset; a variable $G_i$, with the same distribution as $G$, is drawn and this $i$th subset is split into $G_i$ subsets, and so on (see Figure 1).

The $Q$-ary algorithm considered by Mathys and Flajolet [13] corresponds to the case where $G$ is constant and equal to $Q$ and the vector of weights $(V_{i,Q}, 1 \leq i \leq Q)$ is deterministic.

As in Mohamed and Robert [14], for the static case, the key characteristic of this algorithm is a probability distribution $\mathcal{W}$ on $[0,1]$ defined with the branching distribution (the variable $G$), and the weights on each arc [the vector $(V_{1,G}, \ldots, V_{G,G})$]. As we will see, the asymptotic behavior of the algorithm can be described only in terms of the distribution $\mathcal{W}$.

DEFINITION 1. The *splitting measure* is the probability distribution $\mathcal{W}$ on $[0,1]$ defined by, for a nonnegative Borelian function $f$,

$$(10) \quad \int f(x)\mathcal{W}(dx) = \mathbb{E}\left(\sum_{i=1}^{G} V_{i,G} f(V_{i,G})\right) = \sum_{\ell=2}^{+\infty} \sum_{i=1}^{\ell} \mathbb{P}(G = \ell)\mathbb{E}(V_{i,\ell} f(V_{i,\ell})).$$

In the context of fragmentation processes, the measure is related to the *dislocation measure* (see Bertoin [3]).

The following conditions will be assumed throughout the paper:

ASSUMPTIONS (H).

($H_1$) *There exists some* $0 < \delta < 1$ *such that the relation* $\mathcal{W}([0,\delta]) = 1$ *holds;*

$$(H_2) \qquad \int_0^1 \frac{|\log x|}{x}\mathcal{W}(dx) < +\infty.$$

Condition ($H_1$) implies, in particular, the nondegeneracy of the splitting mechanism,

$$\sup_{\ell \geq 2} \sup_{1 \leq i \leq \ell} V_{i,\ell} \leq \delta < 1.$$

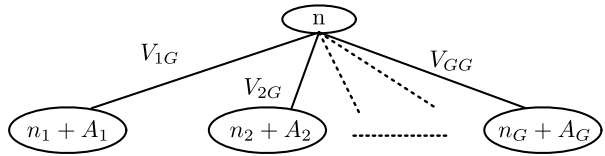

FIG. 1. *First level of the tree algorithm to decompose* $n \geq D$, $n = n_1 + \cdots + n_G$ *and with arrivals* $(A_i)$.

DEFINITION 2.   For $n \geq 0$ and $p \geq 1$, define by $R_n^{A,p}$, the number of nodes with level (generation) less than $p$th when $n$ items are at the root node. The variable $A$ has the same distribution as the common distribution of the i.i.d. sequence $(A_i)$ of the arrivals. By convention, $R_n^{A,\infty} = R_n^A$ and $(R_n^{0,p})$ refers to the static case, that is, when $A \equiv 0$. The variable $G_A$ is defined as

$$G_A = \sum_{i=0}^{+\infty} \sum_{k=1}^{A_i} B_k^i,$$

where, for $i \geq 0$, $(B_k^i, k \geq 0)$ is an i.i.d. Bernoulli sequence independent of $(A_i)$ such that $\mathbb{P}(B_k^i = 1) = \delta^i$.

Note that if $E(A) < +\infty$, the variable $G_A$ is well defined and integrable.

The following lemma establishes a useful relation for tree sizes; the symbol $\leq_{\mathrm{st}}$ is, as usual, for the classical stochastic ordering.

LEMMA 1 (Stochastic inequality).   *For $n \in \mathbb{N}$, under Assumption (H), the relation*

$$(11) \qquad R_n^{A,p} \leq_{\mathrm{st}} R_n^{0,p} + \sum_{i=1}^{R_n^{0,p}} \widetilde{R}_{D-1+G_{A,i}}^{A,p} \mathbf{1}_{\{G_{A,i}>0\}}$$

*holds where $(\widetilde{R}_n^{A,p})$ [resp. $(G_{A,i})$] is a sequence of random variables with the same distribution as $(R_n^{A,p})$ (resp. $G_A$). The variables $(\widetilde{R}_n^{A,p})$, $(R_n^{0,p})$ and $(G_{A,i})$ are independent.*

PROOF.   First note that, for $n \leq m$, one has clearly that $R_n^{A,p} \leq_{\mathrm{st}} R_m^{A,p}$. A coupling is used to prove the stochastic ordering (11). The splitting algorithm is first played only for the initial $n$ items. The total number of nodes up to level $p$ for the associated tree $\mathcal{T}^0$ is $R_n^{0,p}$. The leaves of the tree have, at most, $D$ items, and the internal nodes have more than $D$ items.

Now external arrivals are added to the internal nodes along with the splitting algorithm played on these new items with the branching structure associated with $\mathcal{T}^0$ until they reach one of the leaves of $\mathcal{T}^0$. From there, for all the leaves which have more than $D$ items, the dynamic algorithm is played starting from this node. The number of external items at the leaves has thus to be estimated.

Because of Assumption (H₁), an item in a node containing more than $D$ items is sent to a given child of this node with a probability of, at most, $\delta$. Hence, a given leaf $L = I_i$ of $\mathcal{T}^0$ with depth $i \leq p$ contains at most $D - 1$ initial items and a fraction of the number of new items $A_{L,k}$, $0 \leq k \leq p$, arrives successively at the internal nodes connecting this leaf to the root. Each of the external items arrive at some node and goes to some fixed node below with a probability of, at most, $\delta$, so that, in distribution, there are at most

$$B_1^1 + B_2^1 + \cdots + B_{A_{L,i-1}}^1$$

such items at this node. Similarly, for $2 \leq k \leq i-1$, external items that arrive at a node of level $k$ will reach some fixed node of level $i$ with a probability of, at most, $\delta^{i-k}$. Consequently, the total number $N_L$ of items at leaf $L$ is, in distribution, at most,

$$D - 1 + \sum_{i=1}^{p-2} \sum_{k=1}^{A_{L,i}} B_k^i.$$

If $N_L$ is greater than $D$, the dynamic algorithm continues at that leaf, starting with $N_L + A_0 \leq_{\mathrm{st}} D - 1 + G_A$ items. Note that if this event happens, there must be new items at this leaf, and thus $G_A > 0$. Otherwise, there are only the initial items at $L$, and therefore the algorithm stops.

By noting that the number of leaves of $\mathcal{T}^0$ whose depth $\leq p$ is less than $R_n^{0,p}$, one thus gets the desired relation,

$$R_n^{A,p} \leq_{\mathrm{st}} R_n^{0,p} + \sum_{k=1}^{R_n^{0,p}} \widetilde{R}_{D-1+G_{A,i}}^{A,p} \mathbf{1}_{\{G_{A,i} > 0\}}. \qquad \square$$

THEOREM 3 (Existence of a stable system).  *Under Assumptions $(H)$ for the splitting algorithm, and if $A^\varepsilon$ is a family of integrable integer valued random variables such that*

$$\lim_{\varepsilon \to 0} A^\varepsilon = 0 \ \text{in distribution} \quad \text{and} \quad \limsup_{\varepsilon \to 0} \mathbb{E}(A^\varepsilon \mid A^\varepsilon > 0) < +\infty,$$

*then there exists some $\varepsilon_0 > 0$ and a finite constant $C_{\mathcal{W}}^1$ such that for any $\varepsilon \leq \varepsilon_0$, $R_n^{A^\varepsilon}$ is integrable and*

$$(12) \qquad \mathbb{E}(R_n^{A^\varepsilon}) \leq n C_{\mathcal{W}}^1, \qquad \forall n \geq 0.$$

*In particular, for such an $\varepsilon$, the Poisson transform of the sequence $(\mathbb{E}(R_n^{A^\varepsilon}))$ is defined and differentiable on $\mathbb{R}$.*

Note that the conditions on the family $(A^\varepsilon)$ are quite weak to assert the stability of the algorithm for $\varepsilon$ sufficiently small.

PROOF OF THEOREM 3.  If $F_{A^\varepsilon}$ is a random variable with the same distribution as $(G_{A^\varepsilon} | G_{A^\varepsilon} > 0)$,

$$\mathbb{E}(F_{A^\varepsilon}) = \frac{1}{\mathbb{P}(G_{A^\varepsilon} > 0)} \mathbb{E}(G_{A^\varepsilon} \mathbf{1}_{\{G_{A^\varepsilon} > 0\}})$$

$$\leq \frac{\mathbb{E}(A^\varepsilon)}{\mathbb{P}(A^\varepsilon > 0)(1-\delta)} = \mathbb{E}(A^\varepsilon \mid A^\varepsilon > 0) \frac{1}{1-\delta},$$

the assumptions of the theorem imply that $\mathbb{E}(F_{A^\varepsilon})$ is bounded by some constant $K$ as $\varepsilon$ goes to 0.

By Theorem 3 of Mohamed and Robert [14], there exists a finite constant $C_{\mathcal{W}}^0$ such that

$$\mathbb{E}(R_n^0) \leq n C_{\mathcal{W}}^0, \qquad \forall n \geq 0,$$

and, in particular,

$$\mathbb{E}(R_{D-1+F_{A^\varepsilon}}^0) \leq (D+K)C_{\mathcal{W}}^0.$$

With the variable $R_n^{A^\varepsilon,p}$ being integrable, relation (11) gives the inequality

$$(13) \qquad \mathbb{E}(R_n^{A^\varepsilon,p}) \leq \mathbb{E}(R_n^{0,p}) + \mathbb{E}(R_n^{0,p})\mathbb{P}(G_{A^\varepsilon}>0)\mathbb{E}(R_{D-1+F_{A^\varepsilon}}^{A^\varepsilon,p}), \qquad n \geq 0,$$

and therefore the relation

$$\mathbb{E}(R_{D-1+F_{A^\varepsilon}}^{A^\varepsilon,p}) \leq \mathbb{E}(R_{D-1+F_{A^\varepsilon}}^{0,p})(1+\mathbb{P}(G_{A^\varepsilon}>0)\mathbb{E}(R_{D-1+F_{A^\varepsilon}}^{A^\varepsilon,p}))$$
$$\leq (D+K)C_{\mathcal{W}}^0(1+\mathbb{P}(G_{A^\varepsilon}>0)\mathbb{E}(R_{D-1+F_{A^\varepsilon}}^{A^\varepsilon,p}))$$

holds. It is easy to check that the variables $(G_{A^\varepsilon})$ converge to 0 in distribution as $\varepsilon$ goes to 0. Consequently, there exists some $\varepsilon_0 > 0$ such that if $\varepsilon \leq \varepsilon_0$, then $\mathbb{P}(G_{A^\varepsilon}>0)(D+K)C_{\mathcal{W}}^0 < 1/2$, so that

$$\mathbb{E}(R_{D-1+F_{A^\varepsilon}}^{A^\varepsilon,p}) \leq 2(D+K)C_{\mathcal{W}}^0,$$

and by letting $p$ go to infinity, one gets $\mathbb{E}(R_{D-1+F_{A^\varepsilon}}^{A^\varepsilon}) \leq 2(D+K)C_{\mathcal{W}}^0$. By using again (13) and Theorem 3 of Mohamed and Robert [14], this last inequality gives the relation

$$\mathbb{E}(R_n^{A^\varepsilon,p}) \leq n C_{\mathcal{W}}^0(1+2(D+K)C_{\mathcal{W}}^0),$$

the theorem is proved.  □

Theorem 3 shows that for $\varepsilon$ sufficiently small, the random variables $(R_n^{A^\varepsilon})$ are integrable but also that the variable $R_{A_1^\varepsilon}^{A^\varepsilon}$ is integrable, and, in particular, the Markov chain defined by the transitions (2) is ergodic.

COROLLARY 4 (Stability region for tree algorithms with Poisson arrivals). *When arrivals are Poisson with parameter $\lambda$ and the branching mechanism defined by $\mathcal{W}$ satisfies Assumptions (H), there exists $\lambda_{\mathcal{W}} > 0$ such that:*

(1) *if $\lambda < \lambda_{\mathcal{W}}$, the random variables $(R_n^\lambda)$ are integrable;*
(2) *if $\lambda > \lambda_{\mathcal{W}}$ then $\mathbb{E}(R_n^\lambda) = +\infty$ for all $n \geq D$.*

Proof. If $A_1^\lambda$ is a random variable with a Poisson distribution with parameter $\lambda$, the family $(A_1^\lambda)$ clearly satisfies the assumptions of the above theorem. Hence, there exists $\lambda_0 > 0$ and a constant $C$ such that $\mathbb{E}(R_n^\lambda) < nC$ for all $n \geq 0$.

The sum of the components of the Markov chain defined by the transitions (2) decreases at most of $D-1$ during a time unit and new arrivals have mean $\lambda$, therefore, if $\lambda > D-1$, then $\mathbb{E}(R_n^\lambda) = +\infty$ for all $n \geq D$. The quantity

$$\lambda_{\mathcal{W}} = \sup\{\lambda \geq 0 : \mathbb{E}(R_n^\lambda) < +\infty\}$$

is thus positive and finite. Moreover, it does not depend on $n \geq D$; indeed, if $\mathbb{E}(R_n^\lambda) < +\infty$ and for $m \geq D$, if $m \leq n$, clearly $\mathbb{E}(R_m^\lambda) < \mathbb{E}(R_n^\lambda)$. If $m \geq n$, from (1) one obtains that $R_n^\lambda \geq R_m^\lambda \mathbf{1}_{\{n_1+A_1=m\}}$ so that $R_m^\lambda$ is integrable. Since the function $\lambda \to \mathbb{E}(R_n^\lambda)$ is nondecreasing, one obtains that if $\lambda < \lambda_{\mathcal{W}}$ then, for all $n \geq 0$, the variable $R_n^\lambda$ is integrable. □

3. Poisson transform. From now on and for the rest of the paper, it is assumed that the arrivals are Poisson with parameter $\lambda$ and, as before, one writes $R_n^\lambda$ instead of $R_n^A$. The sequence $\mathcal{N} = (t_n)$, $0 \leq t_1 \leq \cdots \leq t_n \leq \cdots$, is assumed to be a Poisson process with intensity 1 and, for $x \geq 0$, $\mathcal{N}([0,x])$ denotes the number of $t_n$'s in the interval $[0,x]$,

$$\mathcal{N}([0,x]) = \inf\{n : t_{n+1} > x\}.$$

The Poisson transform $\phi_r$ of a sequence $(r_n)$ is given by

$$\phi_r(x) = \sum_{n \geq 0} r_n \frac{x^n}{n!} e^{-x} = \mathbb{E}(r_{\mathcal{N}([0,x])}).$$

Provided that this function is well defined on $\mathbb{R}$, formally,

$$\phi_r'(x) = \sum_{n \geq 0} (r_{n+1} - r_n) \frac{x^n}{n!} e^{-x} = \phi_{\Delta r}(x),$$

where $\Delta r = (r_{n+1} - r_n, n \geq 0)$. In other words, the Poisson transform commutes with the differentiation; the derivative of the Poisson transform of $(r_n)$ is the Poisson transform of the (discrete) derivative of $(r_n)$. The following relation is easily checked by induction, for $n \geq 0$ and $x \geq 0$,

$$(14) \qquad \mathbb{E}(r_{n+\mathcal{N}([0,x])}) = \sum_{k=0}^{n} \binom{n}{k} \phi_r^{(k)}(x).$$

The next proposition establishes an important functional equation for the Poisson transform. For convenience, the Poisson transform of $(\mathbb{E}(R_n^\lambda))$ is denoted by $\phi_\lambda$ instead of $\phi_{R^\lambda}$.

PROPOSITION 5 (Poisson transform). *Provided that $\lambda$ is small enough, the Poisson transform $\phi_\lambda(x)$ of the sequence $(R_n^\lambda)$ satisfies the relation*

$$(15) \qquad \phi_\lambda(x) = \mathbb{E}\left(\frac{1}{W_1}\phi_\lambda(\lambda + W_1 x)\right) + 1 - \phi_C(x),$$

*where $W_1$ is a random variable with distribution $\mathcal{W}$, and $C = (C_m)$ is the sequence defined by $C_m = 0$ for $m \geq D$ and, for $0 \leq m < D$,*

$$C_m = \sum_{k=0}^m \binom{m}{k} \mathbb{E}(W_1^{k-1})\phi_\lambda^{(k)}(\lambda).$$

PROOF. From Theorem 3 and relation (12), one gets that there exists some $\lambda_0$ such that $\phi_\lambda$ is defined on $\mathbb{R}$ when $\lambda < \lambda_0$. By using the splitting property of Poisson processes and by including the boundary cases of (1), one gets the relation

$$
\begin{aligned}
(16) \qquad R_{\mathcal{N}([0,x])}^\lambda &\stackrel{\text{dist}}{=} 1 + \sum_{i=1}^G R_{i,\mathcal{N}([xS_{i-1,G}, xS_{i,G}]) + Z_i}^\lambda \\
&\quad - \mathbf{1}_{\{\mathcal{N}([0,x]) < D\}} \sum_{i=1}^G R_{i,\mathcal{N}([xS_{i-1,G}, xS_{i,G}]) + Z_i}^\lambda,
\end{aligned}
$$

where:

- for $1 \leq i \leq G$, $S_{i,G}$ is the $i$th partial sum of the weights

$$S_{i,G} = V_{1,G} + V_{2,G} + \cdots + V_{i,G},$$

  in particular $S_{G,G} = 1$;
- the variables $R_{i,j}$, $i \geq 1$, $j \geq 0$, are independent and $R_{i,n}$ has the same distribution as $R_n$ for any $n \geq 0$;
- $(Z_i)$ is an i.i.d. sequence of Poisson random variables with parameter $\lambda$.

For $k \geq 0$, the homogeneity properties of Poisson processes give

$$
\begin{aligned}
\mathbb{E}&\left(\mathbf{1}_{\{\mathcal{N}([0,x]) = k\}} \sum_{i=1}^G R_{i,\mathcal{N}([xS_{i-1,G}, xS_{i,G}]) + Z_i}^\lambda\right) \\
&= \mathbb{E}\left(\mathbf{1}_{\{\mathcal{N}([0,x]) = k\}} \sum_{i=1}^G R_{\mathcal{N}([0,xV_{i,G}]) + Z_i}^\lambda\right) \\
&= \mathbb{E}\left(\mathbf{1}_{\{\mathcal{N}([0,x]) = k\}} \frac{R_{\mathcal{N}([0,xW_1]) + Z_1}^\lambda}{W_1}\right) \\
&= \mathbb{E}\left(\frac{R_{\mathcal{N}([0,W_1]) + Z_1}^\lambda}{W_1}\bigg| \mathcal{N}([0,1]) = k\right)\frac{x^k}{k!}e^{-x},
\end{aligned}
$$

where $W_1$ is a random variable whose distribution is $\mathcal{W}$. By taking the expected value of (16), one gets the relation

$$\phi_\lambda(x) = 1 + \mathbb{E}\left(\sum_{i=1}^{G} \phi_\lambda(\lambda + xV_{i,G})\right)$$
$$- \sum_{k=0}^{D-1} \mathbb{E}\left(\frac{R^\lambda_{\mathcal{N}([0,W_1])+Z_1}}{W_1}\bigg|\mathcal{N}([0,1]) = k\right)\frac{x^k}{k!}e^{-x};$$

consequently,

$$\phi_\lambda(x) = 1 + \mathbb{E}\left(\frac{1}{W_1}\phi_\lambda(\lambda + xW_1)\right)$$
$$- \sum_{k=0}^{D-1} \mathbb{E}\left(\frac{R^\lambda_{\mathcal{N}([0,W_1])+Z_1}}{W_1}\bigg|\mathcal{N}([0,1]) = k\right)\frac{x^k}{k!}e^{-x}.$$

For $k \geq 0$, from (14),

$$\mathbb{E}\left(\frac{R^\lambda_{\mathcal{N}([0,W_1])+Z_1}}{W_1}\bigg|\mathcal{N}([0,1]) = k, W_1\right)$$
$$= \sum_{\ell=0}^{k}\binom{k}{\ell}W_1^{\ell-1}(1-W_1)^{k-\ell}\mathbb{E}(R^\lambda_{\ell+Z_1})$$
$$= \sum_{\ell=0}^{k}\sum_{m=0}^{\ell}\binom{k}{\ell}\binom{\ell}{m}W_1^{\ell-1}(1-W_1)^{k-\ell}\phi_\lambda^{(m)}(\lambda)$$
$$= \sum_{m=0}^{k}\binom{k}{m}W_1^{m-1}\phi_\lambda^{(m)}(\lambda).$$

The proposition is proved.  □

*An auto-regressive process with moving average.* At this point, it is natural to introduce the following sequence of random variables.

DEFINITION 6. The process $(X_n^x)$ is defined by $X_0^x = x$, and

(17) $$X_n^x = W_n X_{n-1}^x + 1, \qquad n \geq 1,$$

where $(W_n)$ is an i.i.d. sequence with the same distribution as $W_1$.

The sequence $(X_n^x)$ is an auto-regressive process with moving average. These processes have interesting theoretical properties and play an important role in many areas. In the following the upper index $x$ may be omitted when $x = 0$.

The terms of the sequence $(X_n^x)$ can be expressed as, for $n \geq 0$,

$$X_n^x = x \prod_{i=1}^n W_i + \sum_{p=1}^n \prod_{i=p+1}^n W_i = \pi_n x + X_n^0$$

with, for $1 \leq k$, $\pi_k = \prod_{i=1}^k W_i$. The distribution of the sequence $(W_i)$ being exchangeable, that is, invariant under permutations, one has

$$(18) \qquad X_n \overset{\text{dist}}{=} X_n^* \overset{\text{def}}{=} \sum_{p=0}^{n-1} \pi_p.$$

The sequence $(X_n^*)$ converges almost surely to $X_\infty^*$, and therefore $(X_n^x)$ converges in distribution to the random variable $X_\infty$ such that

$$X_\infty \overset{\text{dist}}{=} W_1 X_\infty + 1 \quad \text{or} \quad X_\infty \overset{\text{dist}}{=} X_\infty^* = \sum_{p=0}^{+\infty} \pi_p$$

the distribution of $X_\infty$ is not explicitly known in general. With this notation, (15) can be rewritten as

$$\phi_\lambda(x) = \mathbb{E}\left( \frac{1}{W_1} \phi_\lambda(\lambda X_1^{x/\lambda}) \right) + 1 - \phi_C(x).$$

By differentiating with respect to $x$, one gets

$$(19) \qquad \phi_\lambda'(x) = \mathbb{E}(\phi_\lambda'(\lambda X_1^{x/\lambda})) - \phi_{\Delta C}(x).$$

Equation (19) expresses $\phi_\lambda'$ as the solution of the Poisson equation associated with the Markov chain $(X_n^{x/\lambda})$ and the function $\phi_{\Delta C}$. Note that, nevertheless, the function $\phi_{\Delta C}$ is depending on $\phi_\lambda$ through its successive derivatives at $\lambda$. By taking $x = X_\infty$ in (19) and integrating, one gets

$$(20) \qquad \mathbb{E}(\phi_{\Delta C}(\lambda X_\infty)) = 0.$$

The iteration of (19) shows that, for $n \geq 1$,

$$\phi_\lambda'(x) = \mathbb{E}(\phi_\lambda'(\lambda X_n^{x/\lambda})) - \sum_{k=0}^{n-1} \mathbb{E}(\phi_{\Delta C}(\pi_k x + \lambda X_k)),$$

and consequently,

$$(21) \qquad \phi_\lambda'(x) = C_\infty - \sum_{k=0}^{+\infty} \mathbb{E}(\phi_{\Delta C}(\pi_k x + \lambda X_k))$$

with $C_\infty$ defined as $\mathbb{E}(\phi'_\lambda(\lambda X_\infty))$. For $k \geq 0$, by using relation (20) and the exchangeability property,

$$\mathbb{E}\left[\frac{1}{\pi_k}(\phi_C(\pi_k x + \lambda X_k) - \phi_C(\lambda X_k))\right]$$

$$= \mathbb{E}\left[\frac{1}{\pi_k}(\phi_C(\pi_k x + \lambda X_k^*) - \phi_C(\lambda X_k^*) - \pi_k x \phi_{\Delta C}(\lambda X_\infty^*))\right],$$

and since $|X_\infty^* - X_k^*| \leq \pi_k/(1-\delta)$ and $X_k^* \leq 1/(1-\delta)$ by assumption (H$_1$),

$$|\phi_C(\pi_k x + \lambda X_k^*) - \phi_C(\lambda X_k^*) - \pi_k x \phi_{\Delta C}(\lambda X_\infty^*)|$$

$$\leq \frac{1}{2}(\pi_k x)^2 \|\phi_{\Delta^2 C}\|_\infty + \pi_k x |\phi_{\Delta C}(\lambda X_\infty^*) - \phi_{\Delta C}(\lambda X_k^*)|$$

$$\leq \pi_k^2\left(\frac{x^2}{2} + \frac{x}{1-\delta}\right)\|\phi_{\Delta^2 C}\|_\infty.$$

The integration of (21), term by term, is therefore valid, and this finally gives the following representation.

PROPOSITION 7 (Representation of Poisson transform). *Provided that $\lambda$ is small enough, the Poisson transform $\phi_\lambda(x)$ of the sequence $(R_n^\lambda)$ satisfies the relation*

$$(22) \quad \phi_\lambda(x) = 1 + x C_\infty + \mathbb{E}\left(\sum_{k=0}^{+\infty} \frac{1}{\pi_k}[\phi_C(\lambda X_k) - \phi_C(\pi_k x + \lambda X_k)]\right),$$

*where $C = (C_n)$ is the sequence defined in Proposition 5 and $C_\infty = \mathbb{E}(\phi'_\lambda(\lambda X_\infty))$;*

(1) *$(W_n)$ is a sequence of i.i.d. random variables whose distribution is $\mathcal{W}$.*
(2) *$(X_n)$ is the auto-regressive process defined by $X_n = W_n X_{n-1} + 1$ for $n \geq 1$ with $X_0 = 0$ and $X_\infty$ is its limit in distribution.*

4. **Stability condition.** In order to get the condition to get the existence of a first moment for the sequence $(R_n^\lambda)$, one has to establish an appropriate representation of this sequence by inverting probabilistically its Poisson transform and to get an expression for the unknown constants $C_0$, $C_1, \ldots, C_{D-1}$ and $C_\infty$.

The notation of Proposition 7 are used. Let $\mathcal{F}_k$ be the $\sigma$-field generated by the random variables $W_1, \ldots, W_k$, and $\mathcal{N}_1$ is another Poisson process with rate 1 independent of $\mathcal{N}$ and $(W_n)$, then, for $k \geq 1$,

$$\phi_C(\pi_k x + \lambda X_k)$$

$$= \mathbb{E}(C_{\mathcal{N}([0,x\pi_k]) + \mathcal{N}_1([0,\lambda X_k])} \mid \mathcal{F}_k)$$

$$= \sum_{m \geq 0} \mathbb{E}(C_{\mathcal{N}([0, x\pi_k]) + \mathcal{N}_1([0, \lambda X X_k])} \mid \mathcal{F}_k, \mathcal{N}([0, x]) = m) \frac{x^m}{m!} e^{-x}$$

$$= \sum_{m \geq 0} \mathbb{E}(C_{\mathcal{N}([0, \pi_k]) + \mathcal{N}_1([0, \lambda X X_k])} \mid \mathcal{F}_k, \mathcal{N}([0, 1]) = m) \frac{x^m}{m!} e^{-x}.$$

With (22), one obtains the relation

$$\phi_\lambda(x) = 1 + x C_\infty$$
$$+ \sum_{m \geq 0} \frac{x^m}{m!} e^{-x} \sum_{k \geq 0} \mathbb{E}\bigg( \frac{1}{\pi_k} [C_{\mathcal{N}([0, \lambda X X_k])}$$
$$- C_{\mathcal{U}_m([0, \pi_k]) + \mathcal{N}([0, \lambda X X_k])}] \bigg),$$

where, if $U_1, \ldots, U_n$ are $n$ i.i.d. random variables uniformly distributed on $[0, 1]$, for $0 \leq x \leq 1$, $\mathcal{U}_n([0, x])$ denotes the number of $U_k$'s in the interval $[0, x]$. These variables are ordered as $U_{(1)}^n \leq U_{(2)}^n \leq \cdots \leq U_{(n)}^n$, in particular, for $m \geq 1$, $\{\mathcal{U}_n([0, x]) \geq m\} = \{U_{(m)}^n \leq x\}$. By identifying the coefficients of the above expression, one gets the following proposition.

PROPOSITION 8 (Representation of the average size of the tree).   *Under Assumptions (H) and for $\lambda$ sufficiently small,*

(23)
$$\mathbb{E}(R_n^\lambda) = 1 + n C_\infty$$
$$+ \sum_{k \geq 0} \mathbb{E}\bigg( \frac{1}{\pi_k} [C_{\mathcal{N}([0, \lambda X X_k])} - C_{\mathcal{U}_n([0, \pi_k]) + \mathcal{N}([0, \lambda X X_k])}] \bigg), \qquad n \geq 0,$$

*where $\mathcal{N}$ is a Poisson process with rate 1 and:*

(1) *$C = (C_n)$ is the sequence defined in Proposition 5 and $C_\infty = \mathbb{E}(\phi_\lambda'(\lambda X_\infty))$.*
(2) *For uniformly distributed random variables $(U_i, 1 \leq i \leq n)$ on $[0, 1]$ and $0 \leq x \leq 1$, the quantity $\mathcal{U}_n([0, x])$ denotes the number of $U_i$'s in the interval $[0, x]$.*

*Determination of the constants.*   In order to get an explicit representation of the sequence $(\mathbb{E}(R_n^\lambda))$, the $D$ unknown coefficients $C_0, \ldots, C_{D-1}$ (recall that the other $C_k$'s are null) and the constant $C_\infty = \mathbb{E}(\phi_\lambda'(\lambda X_\infty))$ have to be determined. The method used by Fayolle, Flajolet and Hofri [8] to determine these coefficients in the binary case apparently cannot be extended to other tree structures.

(i) The boundary conditions $\mathbb{E}(R_m^\lambda) = 1$ for $1 \le m \le D-1$ translate into $D-1$ linear equations involving these $D+1$ unknown constants,

$$(24) \qquad \sum_{\ell=0}^{D-1} M_{m,\ell}^\lambda C_\ell + M_{m,D}^\lambda C_\infty = 0, \qquad 1 \le m \le D-1,$$

with, for $1 \le m \le D-1, 0 \le \ell \le D-1$,

$$M_{m,\ell}^\lambda = \sum_{k \ge 0} \mathbb{E}\left( \frac{1}{\pi_k} [\mathbf{1}_{\{\mathcal{N}([0,\lambda X_k]) = \ell\}} - \mathbf{1}_{\{\mathcal{U}_m([0,\pi_k]) + \mathcal{N}([0,\lambda X_k]) = \ell\}}] \right)$$

and $M_{m,D}^\lambda = m$.

(ii) Equation (20) gives the additional relation

$$(25) \qquad M_{D,0}^\lambda C_0 + M_{D,1}^\lambda C_1 + \cdots + M_{D,D-1}^\lambda C_{D-1} = 0$$

with

$$M_{D,\ell}^\lambda = \mathbb{E}\left( \frac{(\lambda X_\infty)^\ell}{\ell!} (\lambda^{-1} \ell X_\infty^{-1} - 1) e^{-\lambda X_\infty} \right), \qquad 0 \le \ell \le D-1,$$

and $M_{D,D}^\lambda = 0$.

(iii) The final equation is obtained by plugging $x = \lambda$ in (22) so that, since $C_0 = \mathbb{E}(G)\phi(\lambda)$ by Proposition 5,

$$
\begin{aligned}
-1 = {}& \lambda C_\infty - \frac{1}{\mathbb{E}(G)} C_0 \\
(26) \quad & + \sum_{m=0}^{D-1} C_m \mathbb{E}\left( \sum_{k=0}^{+\infty} \frac{1}{\pi_k} \left[ \frac{(\lambda X_k)^m}{m!} e^{-(\lambda X_k)} \right. \right. \\
& \qquad\qquad \left. \left. - \frac{(\pi_k \lambda + \lambda X_k)^m}{m!} e^{-(\pi_k \lambda + \lambda X_k)} \right] \right).
\end{aligned}
$$

*The matrix $M_\lambda$.* The square matrix $M_\lambda = (M_{m,\ell}^\lambda, 1 \le m \le D+1, 0 \le \ell \le D)$ is defined as follows:

$$
\begin{cases}
M_{m,\ell}^\lambda = \sum_{k \ge 0} \mathbb{E}\left( \frac{1}{\pi_k} [\mathbf{1}_{\{\mathcal{N}([0,\lambda X_k])=\ell\}} - \mathbf{1}_{\{\mathcal{U}_m([0,\pi_k])+\mathcal{N}([0,\lambda X_k])=\ell\}}] \right), \\
\quad m < D, \ell \ne D, \\
M_{D,\ell}^\lambda = \mathbb{E}\left( \frac{(\lambda X_\infty)^\ell}{\ell!} [\lambda^{-1} \ell X_\infty^{-1} - 1] e^{-\lambda X_\infty} \right), \\
\quad 0 \le \ell \le D-1, \\
M_{D+1,\ell}^\lambda = \mathbb{E}\left( \sum_{k=0}^{+\infty} \frac{1}{\pi_k} [X_k^\ell - (\pi_k + X_k)^\ell e^{-\lambda \pi_k}] \frac{\lambda^\ell}{\ell!} e^{-\lambda X_k} \right), \\
\quad 1 \le \ell \le D-1, \\
M_{D+1,0}^\lambda = \mathbb{E}\left( \sum_{k=0}^{+\infty} \frac{1}{\pi_k} [1 - e^{-\lambda \pi_k}] e^{-\lambda X_k} \right) - \frac{1}{\mathbb{E}(G)}, \\
M_{D,D}^\lambda = 0, \qquad M_{D+1,D}^\lambda = \lambda.
\end{cases}
$$

By gathering equations (24), (25) and (26) and denoting $e_{D+1} = (0,0,\ldots,0,1)$ and by $\underline{C} = (C_0, C_1, \ldots, C_{D-1}, C_\infty)$, the vector of constants, one gets the linear relation

(27)                        $$M_\lambda \cdot C = -e_{D+1}.$$

The following theorem is the main result concerning the ergodicity of the tree algorithm.

THEOREM 9 (Stability of tree algorithm with Poisson arrivals). *Under Assumptions (H) for the splitting distribution $\mathcal{W}$, if $M_\lambda$ is the matrix defined above and*

$$\lambda_c = \inf\{\lambda > 0 : \det M_\lambda = 0\},$$

*then $\lambda_c > 0$ and for any $\lambda < \lambda_c$, the size of the tree associated to the tree algorithm is integrable.*

PROOF. With the same notation as before, for $\lambda = 0$, then $X_k = 0$ for $0 \le k \le +\infty$ and the $D$th and $(D+1)$th rows of the matrix $\mathcal{M}_\lambda$ are given by

$$M_D = (-1, 1, 0, \ldots, 0) \quad \text{and} \quad M_{D+1} = (-1/\mathbb{E}(G), 0, \ldots, 0, 0)$$

by expanding, with respect to these rows, one gets

$$
\det M_0 = \frac{1}{\mathbb{E}(G)} \begin{vmatrix}
M_{1,2}^0 & M_{1,3}^0 & \cdots & M_{1,D-1}^0 & 1 \\
M_{2,2}^0 & M_{2,3}^0 & \cdots & M_{2,D-1}^0 & 2 \\
\vdots & \vdots & \vdots & \vdots & \vdots \\
M_{D-1,2}^0 & M_{D-1,3}^0 & \cdots & M_{D-1,D-1}^0 & D-1
\end{vmatrix}.
$$

Since, for $1 \leq m, \ell \leq D - 1$,

$$M_{m,\ell}^0 = - \sum_{k \geq 0} \mathbb{E}\left( \frac{1}{\pi_k} \mathbf{1}_{\{ \mathcal{U}_m([0, \pi_k]) = \ell \}} \right);$$

then, $M_{m,\ell}^0 = 0$ for $\ell > m$, and hence

$$\det \mathcal{M}_0 = \frac{1}{\mathbb{E}(G)} M_{2,2}^0 \cdots M_{D-1,D-1}^0 \neq 0.$$

Due to the explicit expression of the matrix $M_\lambda$, the function $\lambda \to \det M_\lambda$ is clearly continuous so that $\lambda_c > 0$.

Corollary 4 shows the existence of some constant $\lambda_{\mathcal{W}}$ such that, for $\lambda < \lambda_{\mathcal{W}}$, random variables $(R_n^\lambda)$ are integrable and their expected values are given by (23) and the constant vector $C$ in this expression satisfies (27). Hence, for $\lambda < \lambda_{\mathcal{W}} \wedge \lambda_c$, there exists a unique $C$ such that (23) holds for $\mathbb{E}(R_n^\lambda)$ for $n \geq 0$. Since the function $\lambda \to \mathbb{E}(R_n^\lambda)$ is nondecreasing and because of the existence of a solution to (27) $\lambda < \lambda_c$, the expression given by (23) is finite for any $\lambda < \lambda_c$; thus one concludes necessarily, by Corollary 4, that $\lambda_c \leq \lambda_{\mathcal{W}}$. The theorem is proved. □

REMARKS.

(1) It is very likely that $\lambda_{\mathcal{W}}$ defined in Proposition 4 and $\lambda_c$ are equal, but we have not been able to prove it. For $\lambda = \lambda_c$, at least one of the determinants of the Cramér formula should be nonzero which would imply that at least one of the $(C_k)$'s is infinite, and therefore that the random variables $(R_n^\lambda)$ are not integrable.

(2) The coefficients of the matrix $M_\lambda$ are expressed in terms of the distribution of the auto-regressive process $(X_n)$. An explicit, usable representation of this distribution is available mostly through Laplace transform functionals, the invariant distribution included.

  Although it is not easy to handle, the Introduction of the auto-regressive process is, in our opinion, the key ingredient in our analysis. It plays a major role in representation (23) of the sequence $(\mathbb{E}(R_n^\lambda))$. It should also be kept in mind that one relation used to determine the constants $(C_k)$ is (20) which comes directly from the fact that $(X_n)$ has an equilibrium distribution. In a purely analytic setting (i.e., without this probabilistic representation), an analogous equation would probably require some spectral analysis in a functional space.

EXAMPLES.

(a) *Static case.* In this case $\lambda = 0$, and the components of the vector $C = (C_i, 0 \leq i \leq D - 1)$ are constant and equal to the average branching degree $\mathbb{E}(G) = \mathbb{E}(W_1^{-1})$. For $n \geq D$,

$$\mathbb{E}(R_n^0) = 1 + \mathbb{E}(G) \sum_{k \geq 0} \mathbb{E}\left(\frac{1}{\pi_k}\mathbf{1}_{\{U_{(D)}^n \leq \pi_k\}}\right),$$

which is the expression established for the static splitting algorithm in Mohamed and Robert [14].

(b) *Binary tree algorithm.* In the binary case, $G \equiv 2$, $D = 2$ and the splitting measure is $\mathcal{W} = p\delta_p + q\delta_q$. In this case the average cost of the algorithm is expressed as follows:

$$\mathbb{E}(R_n^\lambda) = 1 + nC_\infty - \phi_\lambda'(\lambda) \sum_{k \geq 0} \mathbb{E}\left(\frac{e^{-\lambda X_k}}{\pi_k}\mathbf{1}_{\{U_{(1)}^n \leq \pi_k\}}\right)$$

$$+ (2\phi_\lambda(\lambda) + \phi_\lambda'(\lambda))\left(\sum_{k \geq 0} \mathbb{E}\left(\frac{\lambda X_k e^{-\lambda X_k}}{\pi_k}\mathbf{1}_{\{U_{(1)}^n \leq \pi_k\}}\right)\right.$$

$$\left. + \sum_{k \geq 0} \mathbb{E}\left(\frac{e^{-\lambda X_k}}{\pi_k}\mathbf{1}_{\{U_{(2)}^n \leq \pi_k\}}\right)\right).$$

The two coefficients $C_0$ and $C_1$ satisfy

$$C_0 = 2\phi_\lambda(\lambda), \qquad C_1 = 2\phi_\lambda(\lambda) + \phi_\lambda'(\lambda).$$

Equation (20) implies that

$$[\mathbb{E}(e^{-\lambda X_\infty}) - \lambda\mathbb{E}(X_\infty e^{-\lambda X_\infty})]C_1 = \mathbb{E}(e^{-\lambda X_\infty})C_0,$$

which gives the relation

$$(28) \qquad \phi_\lambda'(\lambda) = 2\left(\frac{\mathbb{E}(e^{-\lambda X_\infty})}{\mathbb{E}(e^{-\lambda X_\infty}) - \lambda\mathbb{E}(X_\infty e^{-\lambda X_\infty})} - 1\right)\phi_\lambda(\lambda).$$

Note that in the case of the symmetric binary algorithm, the limit $X_\infty$ is constant and equal to 2.

A identity similar to (28) has been established in Fayolle et al. [9]. By taking advantage of the fact that if one plugs $x = \lambda/p$ and $x = \lambda/q$ successively into (15) [(6) in this case], one gets the relation

$$\phi_\lambda'(\lambda) = 2(K-1)\phi_\lambda(\lambda),$$

where

$$K = \begin{cases} -\dfrac{e^{-\lambda/p} - e^{-\lambda/q}}{\lambda/p e^{-\lambda/p} - \lambda/q e^{-\lambda/q}}, & \text{if } p \neq 1/2, \\ 1/(1-2\lambda), & \text{otherwise.} \end{cases}$$

This trick turns out to be specific to binary trees and does not seem to have a generalization for other random trees. Interestingly, when $p \neq 1/2$, the representation of the constant $K$ by (28) gives the following relation for $g(\lambda)$, the Laplace transform of $X_\infty$ at $\lambda$,

$$-\frac{g'(\lambda)}{g(\lambda)} = \frac{(1+\lambda/p)e^{-\lambda/p} - (1+\lambda/q)e^{-\lambda/q}}{\lambda(e^{-\lambda/p} - e^{-\lambda/q})},$$

which can be solved as

$$(29) \qquad \mathbb{E}(e^{-\lambda X_\infty}) = \frac{1}{1/(1-p)-1/p} \frac{e^{-\lambda/p} - e^{-\lambda/(1-p)}}{\lambda},$$

which gives an explicit representation of the Laplace transform of the invariant measure of the auto-regressive process in this case.

**5. Asymptotic analysis of the average size of the tree.** In this section, the asymptotic behavior of the sequence $(\mathbb{E}(R_n^\lambda))$ is investigated for $\lambda < \lambda_c$ where $\lambda_c$ is defined in Theorem 9. The goal is to establish an analogue of the law of large numbers for these expected values. As noted before, Fayolle, Flajolet and Hofri [8] (for the binary tree) is the only rigorous result we know in this domain.

Equation (23) gives the representation, for $n \geq D$,

$$(30) \quad \mathbb{E}(R_n^\lambda) = 1 + nC_\infty - \sum_{i=0}^{D-1} \sum_{k \geq 0} \mathbb{E}\left(\frac{1}{\pi_k} \Delta C_{i+\mathcal{N}([0,\lambda X_k])} \mathbf{1}_{\{U_{(i+1)}^n \leq \pi_k\}}\right),$$

where $(\Delta C_i)$ is the sequence $(C_{i+1} - C_i)$ and, for $0 \leq i \leq D-1$ and $n \geq 1$. Recall that $U_{(i)}^n$ is the $i$th smallest term of $n$ independent uniform random variables on $[0,1]$.

In a first step, it is shown that the series associated to $i = 0$ in (30) is vanishing for the asymptotic behavior of the sequence $(\mathbb{E}(R_n^\lambda)/n)$. This is a crucial result since the arguments to derive a law of large numbers rely on an integrability property which is not satisfied for this series.

LEMMA 2. *Under Assumptions (H), the relation*

$$\lim_{n \to +\infty} \frac{1}{n} \sum_{k \geq 0} \mathbb{E}\left(\frac{1}{\pi_k} \Delta C_{\mathcal{N}([0,\lambda X_k])} \mathbf{1}_{\{U_{(1)}^n \leq \pi_k\}}\right) = 0$$

*holds.*

PROOF. Equation (19) gives the relation

$$A_n \stackrel{\text{def}}{=} \sum_{k \geq 0} \mathbb{E}\left(\frac{1}{\pi_k} \Delta C_{\mathcal{N}([0,\lambda X_k])} \mathbf{1}_{\{U_{(1)}^n \leq \pi_k\}}\right)$$

$$= \sum_{k \geq 0} \mathbb{E}\left( \frac{1}{\pi_k} \phi_{\Delta C}(\lambda X_k) \mathbf{1}_{\{U_{(1)}^n \leq \pi_k\}} \right)$$

$$= \sum_{k \geq 0} \mathbb{E}\left( \frac{1}{\pi_k} [\Delta R_{\mathcal{N}([0, \lambda X_{k+1}])}^\lambda - \Delta R_{\mathcal{N}([0, \lambda X_k])}^\lambda] \mathbf{1}_{\{U_{(1)}^n \leq \pi_k\}} \right)$$

$$\leq \sup_{0 \leq x \leq \lambda/(1-\delta)} \phi_{\Delta R^\lambda}(x) \cdot \sum_{k \geq 0} \mathbb{E}\left( \frac{1}{\pi_k} \mathbf{1}_{\{\mathcal{N}([\lambda X_k, \lambda X_{k+1}]) \neq 0, U_{(1)}^n \leq \pi_k\}} \right),$$

by assumption (H$_1$). For $k \geq 0$, by exchangeability of the sequence $(W_i)$ and definition (18), one gets

$$\mathbb{E}\left( \frac{1}{\pi_k} \mathbf{1}_{\{\mathcal{N}([\lambda X_k, \lambda X_{k+1}]) \neq 0, U_{(1)}^n \leq \pi_k\}} \right)$$

$$= \mathbb{E}\left( \frac{\pi_1}{\pi_{k+1}} \mathbf{1}_{\{\mathcal{N}([\lambda X_k^*, \lambda X_{k+1}^*]) \neq 0, U_{(1)}^n \leq \pi_{k+1}/\pi_1\}} \right)$$

$$\leq \delta \mathbb{E}\left( \frac{1}{W_{k+1}} \mathbf{1}_{\{W_1 U_{(1)}^n \leq \delta^{k+1}\}} \right)$$

$$= \delta \mathbb{E}(G) \mathbb{P}(W_1 U_{(1)}^n \leq \delta^{k+1}).$$

By summing up these terms, this gives the following upper bound for $A_n$, for some constant $C$,

$$A_n \leq C \sum_{k \geq 0} \mathbb{P}(W_1 U_{(1)}^n \leq \delta^{k+1}) \leq C \mathbb{E}(\lceil -\log_{1/\delta}(W_1 U_{(1)}^n) \rceil).$$

Since this term is of the order of $\log n$, the sequence $(A_n/n)$ converges to 0. $\square$

Before analyzing the asymptotic behavior of $(\mathbb{E}(R_n^\lambda)/n)$, Propositions 9 and 11 from Mohamed and Robert [14] obtained in the static case are summarized in the following proposition.

PROPOSITION 10.    *Under Assumption (H), for* $i \geq 1$, *if*

$$E_{i,n} \overset{\text{def}}{=} \sum_{k \geq 0} \mathbb{E}\left( \frac{1}{\pi_k} \mathbf{1}_{\{U_{(i+1)}^n \leq \pi_k\}} \right).$$

(1) *If the random variable* $-\log W_1$ *is nonarithmetic, then*

$$\lim_{n \to +\infty} \frac{E_{i,n}}{n} = \frac{\mathbb{E}(G)}{i \mathbb{E}(|\log W_1|)}.$$

(2) *If the random variable* $-\log W_1$ *is arithmetic and* $\xi > 0$ *is its span, then, as* $n$ *goes to infinity,*

$$\lim_{n \to +\infty} \frac{E_{i,n}}{n} - F_i\left(\frac{\log n}{\xi}\right) = 0,$$

*where* $F_i$ *is the periodic function defined by*

(31) $$F_i(x) = \frac{\mathbb{E}(G)}{\mathbb{E}(|\log(W_1)|)} \frac{\xi}{1 - e^{-\xi}} \int_0^{+\infty} \exp\left(-\xi\left\{x - \frac{\log y}{\xi}\right\}\right) \frac{y^{i-1}}{i!} e^{-y} \, dy,$$

*and* $\{x\} = x - \lfloor x \rfloor$.

The next proposition "decouples" the process $(X_n)$ and the counting process associated to the sequence $(\pi_k)$.

PROPOSITION 11. *Under Assumption* (H), *for* $1 \le i \le D - 1$,

$$\lim_{n \to +\infty} \frac{1}{n} \sum_{k \ge 0} \mathbb{E}\left(\frac{1}{\pi_k} \mathbf{1}_{\{U^n_{(i+1)} \le \pi_k\}} [\Delta C_{i + \mathcal{N}([0, \lambda X_k])} - \mathbb{E}(\Delta C_{i + \mathcal{N}([0, \lambda X_\infty])})]\right) = 0.$$

PROOF. By using definition (18) and the exchangeability property, one has, for $p \ge 1$,

$$\frac{1}{n} \sum_{k \ge p} \mathbb{E}\left(\frac{1}{\pi_k} |\Delta C_{i + \mathcal{N}([0, \lambda X_k])} - \Delta C_{i + \mathcal{N}([0, \lambda X_p])}| \mathbf{1}_{\{U^n_{(i+1)} \le \pi_k\}}\right)$$

$$= \frac{1}{n} \sum_{k \ge p} \mathbb{E}\left(\frac{1}{\pi_k} |\Delta C_{i + \mathcal{N}([0, \lambda X_k^*])} - \Delta C_{i + \mathcal{N}([0, \lambda X_p^*])}| \mathbf{1}_{\{U^n_{(i+1)} \le \pi_k\}}\right)$$

$$\le \frac{1}{n} \|\Delta C\|_\infty \sum_{k \ge 0} \mathbb{E}\left(\frac{1}{\pi_k} \mathbb{P}(\mathcal{N}([\lambda X_k^*, \lambda X_p^*]) \ne 0 \mid \mathcal{F}_k) \mathbf{1}_{\{U^n_{(i+1)} \le \pi_k\}}\right)$$

$$= \frac{1}{n} \|\Delta C\|_\infty \mathbb{E}\left(\sum_{k \ge 0} \frac{1}{\pi_k} (1 - e^{-\lambda |X_k^* - X_p^*|}) \mathbf{1}_{\{U^n_{(i+1)} \le \pi_k\}}\right)$$

$$\le \|\Delta C\|_\infty (1 - \exp(-\lambda \delta^p / (1 - \delta))) E_{i,n},$$

by assumption (H$_1$) where $E_{i,n}$ is defined in Proposition 10. One can therefore choose a $p$ sufficiently large so that the above difference is arbitrarily small, uniformly in $n \ge 1$.

One has thus to investigate the asymptotic behavior of

$$\frac{1}{n} \mathbb{E}\left(\Delta C_{i + \mathcal{N}([0, \lambda X_p])} \sum_{k \ge p} \frac{1}{\pi_k} \mathbf{1}_{\{U^n_{(i+1)} \le \pi_k\}}\right) = \mathbb{E}(\Delta C_{i + \mathcal{N}([0, \lambda X_p])} G_p(n))$$

with

$$G_p(n) = \mathbb{E}\left(\frac{1}{n}\sum_{k \geq p}\frac{1}{\pi_k}\mathbf{1}_{\{U_{(i+1)}^n \leq \pi_k\}} \mid \mathcal{F}_p\right).$$

When $-\log W_1$ is nonarithmetic, Proposition 10 shows that $\mathbb{E}(G_p(n))$ converges. With the same argument as in Mohamed and Robert, the conditioning being on the first $p$ elements of the sequence $(W_k)$, almost surely, $G_p(n)$ converges to the same limit as the sequence $(\mathbb{E}(G_p(n)))$. Consequently, by denoting $x^+$, the nonnegative part of $x \in \mathbb{R}$,

$$\mathbb{E}(|\mathbb{E}(G_p(n)) - G_p(n)|) = 2\mathbb{E}([\mathbb{E}(G_p(n)) - G_p(n)]^+),$$

and since the quantity $[\mathbb{E}(G_p(n)) - G_p(n)]^+$ is bounded by $\sup_n \mathbb{E}(G_p(n))$, Lebesgue's theorem gives that the sequence $(\mathbb{E}(G_p(n)) - G_p(n))$ converges to 0 in $L_1$. Therefore, the quantity

$$|\mathbb{E}(\Delta C_{i+\mathcal{N}([0,\lambda X_p])}[G_p(n) - \mathbb{E}[G_p(n)]])| \leq \|\Delta C\|_\infty \mathbb{E}(|G_p(n) - \mathbb{E}(G_p(n))|)$$

converges to 0 as $n$ goes to infinity. The proposition is therefore proved in this case. When $-\log W_1$ is arithmetic with range $\xi\mathbb{N}$, the argument is similar by using the fact that the convergence to 0 of $(G_p(n) - F_i(\log n/\xi))$ holds almost surely and for the expected value.  $\square$

The main result of this section can now be stated. It is a direct consequence of representation (30), Lemma 2, Proposition 10 and Proposition 11.

THEOREM 12.  *If $\lambda < \lambda_c$ defined in Theorem 9 and under Assumption (H),*

(1) *if the random variable $-\log W_1$ is nonarithmetic, then*

$$\lim_{n \to +\infty}\frac{\mathbb{E}(R_n^\lambda)}{n} = C_\infty - \frac{\mathbb{E}(G)}{\mathbb{E}(|\log W_1|)}\sum_{i=1}^{D-1}\frac{1}{i}\mathbb{E}(\Delta C_{i+\mathcal{N}([0,\lambda X_\infty])}).$$

(2) *If the random variable $-\log W_1$ is arithmetic and $\xi > 0$ is its span, then*

$$\lim_{n \to +\infty}\frac{\mathbb{E}(R_n^\lambda)}{n} - C_\infty - \sum_{i=1}^{D-1}\frac{1}{i}\mathbb{E}(\Delta C_{i+\mathcal{N}([0,\lambda X_\infty])})F_i\left(\frac{\log n}{\xi}\right) = 0$$

$(F_i, 1 \leq i \leq D-1)$ *being the periodic functions defined by (31)*

*where:*

– *for $i \geq 1$, $\Delta C_i = C_{i+1} - C_i$ with $\underline{C} = (C_0, C_1, \ldots, C_{D-1}, C_\infty)$ being the vector solution of the equation*

$$M_\lambda \cdot \underline{C} = -e_{D+1}$$

*with $C_k = 0$ for $k \geq D$ and $M_\lambda$ is the matrix above equation (27).*

– $\mathcal{N}$ *is a Poisson point process with rate* 1.
– *The variable* $X_\infty$ *has the invariant distribution of the auto-regressive process* $(X_n)$ *defined by* $X_{n+1} = W_n X_n + 1$, $n \geq 0$.

LMV
Universite de Versailles-saint-Quentin-en-Yvelines
78035 Versailles Cedex
France
E-mail: Hanene.Mohamed@math.uvsq.fr
URL: http://www.math.uvsq.fr/~mohamed/

INRIA
domaine de Voluceau, B.P. 105
78153 Le Chesnay Cedex
France
E-mail: Philippe.Robert@inria.fr
URL: http://www-rocq.inria.fr/~robert